\documentclass{gtart}
\gtart

\usepackage{amsthm}
\usepackage{amsgen}
\usepackage{amsmath,amssymb}
\usepackage{curvesls}
\usepackage{epic}
\usepackage[all,color]{xy}
\usepackage[dvips]{graphicx}
\usepackage{psfrag}
\usepackage[T1]{fontenc}
\usepackage[sc]{mathpazo}
\usepackage{color, xcolor, colortbl}
\usepackage{float}
\usepackage[colorlinks=true, urlcolor=blue, linkcolor=red]{hyperref}
\usepackage{arydshln, multirow} %% dashed lines in tables
\usepackage{multicol}
\usepackage{fancybox}
\newtheorem{thm}{Theorem}[section]
\newtheorem{lem}[thm]{Lemma}
\newtheorem{cor}[thm]{Corollary}
\newtheorem{prop}[thm]{Proposition}

\newtheorem{conj}[thm]{Conjecture}

\theoremstyle{definition}
\newtheorem{defn}[thm]{Definition}

\newcommand{\Real}{{\mathbb R }}
\newcommand{\Rat}{{\mathbb Q}}

\newcommand{\Zed}{{\mathbb Z }}

\newcommand{\Diff}{{\mathrm{Diff}}}
\newcommand{\Map}{{\mathrm{Map}}}

\newcommand{\Emb}{{\mathrm{Emb}}}

\newcommand\scalemath[2]{\scalebox{#1}{\mbox{\ensuremath{\displaystyle #2}}}}

\makeatletter
\renewcommand*\env@matrix[1][\arraystretch]{%
  \edef\arraystretch{#1}%
  \hskip -\arraycolsep
  \let\@ifnextchar\new@ifnextchar
  \array{*\c@MaxMatrixCols c}}
\makeatother

\newcommand\cone[1]{ \scalemath{0.6}{\begin{matrix}[0.5]#1 \end{matrix}} }
\newcommand\ctwo[2]{ \scalemath{0.6}{\begin{matrix}[0.5]#1 \\ #2\end{matrix}} }
\newcommand\cthree[3]{\scalemath{0.6}{\begin{matrix}[0.5]#1 \\ #2 \\ #3\end{matrix}}}
\newcommand\cfour[4]{\scalemath{0.6}{\begin{matrix}[0.5]#1 \\ #2 \\ #3 \\ #4\end{matrix}}}

\begin{document}

\title{A $2$-torsion invariant of $2$-knots}
\author{Ryan Budney}

\address{
Mathematics and Statistics, University of Victoria\\
PO BOX 3045 STN CSC, Victoria, B.C., Canada V8W 3P4
}
\email{rybu@uvic.ca}

\begin{abstract} 
In this paper we describe what should perhaps be called a `type-2' Vassiliev invariant of knots $S^2 \to S^4$.
We give a formula for an invariant of $2$-knots, taking values in $\Zed_2$ that can be computed 
in terms of the double-point diagram of the knot.  The double-point diagram is a collection of curves
and diffeomorphisms of curves, in the domain $S^2$, that describe the crossing data with respect to a projection, 
analogous to a chord diagram for a projection of a classical knot $S^1 \to S^3$.  Our formula turns the computation of the 
invariant into a planar geometry problem. 
More generally, we describe a numerical invariant of families of knots $S^j \to S^n$, 
for all $n \geq j+2$ and $j \geq 1$. In the co-dimension two case $n=j+2$ the invariant is an isotopy invariant,
and either takes values in $\Zed$ or $\Zed_2$ depending on a parity issue.  
\end{abstract}

\maketitle

\section{Introduction}\label{intro}

In this paper we describe a $2$-torsion isotopy invariant of $2$-knots, i.e. smooth embeddings $S^2 \to S^4$. 
We conjecture this is not a new invariant.  Our interest in this invariant comes from its form, specifically how it is
computed, involving the geometry of circles.  This connection makes our invariant analogous to the type-2 invariant of
classical knots \cite{BCSS}.

The invariant will be defined in a language of {\it configuration spaces} $C_k(M)$ and
the geometry of circles. See the Definition \ref{condef} for how we use these terms.
Given a $2$-knot $f : S^2 \to S^4$ the submanifold $\mathcal C_f \subset C_5(S^2)$ is
defined by the condition that $p \in C_5(S^2)$ belongs to $\mathcal C_f$ if and only if the 
points $p = (p_1, \cdots, p_5)$ and $f_*(p) = (f(p_1), \cdots, f(p_5)) \in C_5(S^4)$
satisfy the four conditions below.  
\begin{enumerate}
\item The points $p \in C_5(S^2)$ sit on a round circle in $S^2$.
\item The points $f_*(p) \in C_5(S^4)$ sit on a round circle in $S^4$
\item The points $p$ are in the cyclic order of the circle they lie on.
\item The points $f_*(p)$ are in non-consecutive order in the circle they lie on.
\end{enumerate}

\begin{figure}[H]
{
\psfrag{p1}[tl][tl][0.7][0]{$p_1$}
\psfrag{p2}[tl][tl][0.7][0]{$p_2$}
\psfrag{p3}[tl][tl][0.7][0]{$p_3$}
\psfrag{p4}[tl][tl][0.7][0]{$p_4$}
\psfrag{p5}[tl][tl][0.7][0]{$p_5$}
\psfrag{f1}[tl][tl][0.7][0]{$f(p_1)$}
\psfrag{f2}[tl][tl][0.7][0]{$f(p_3)$}
\psfrag{f3}[tl][tl][0.7][0]{$f(p_5)$}
\psfrag{f4}[tl][tl][0.7][0]{$f(p_2)$}
\psfrag{f5}[tl][tl][0.7][0]{$f(p_4)$}
$$\includegraphics[width=9cm]{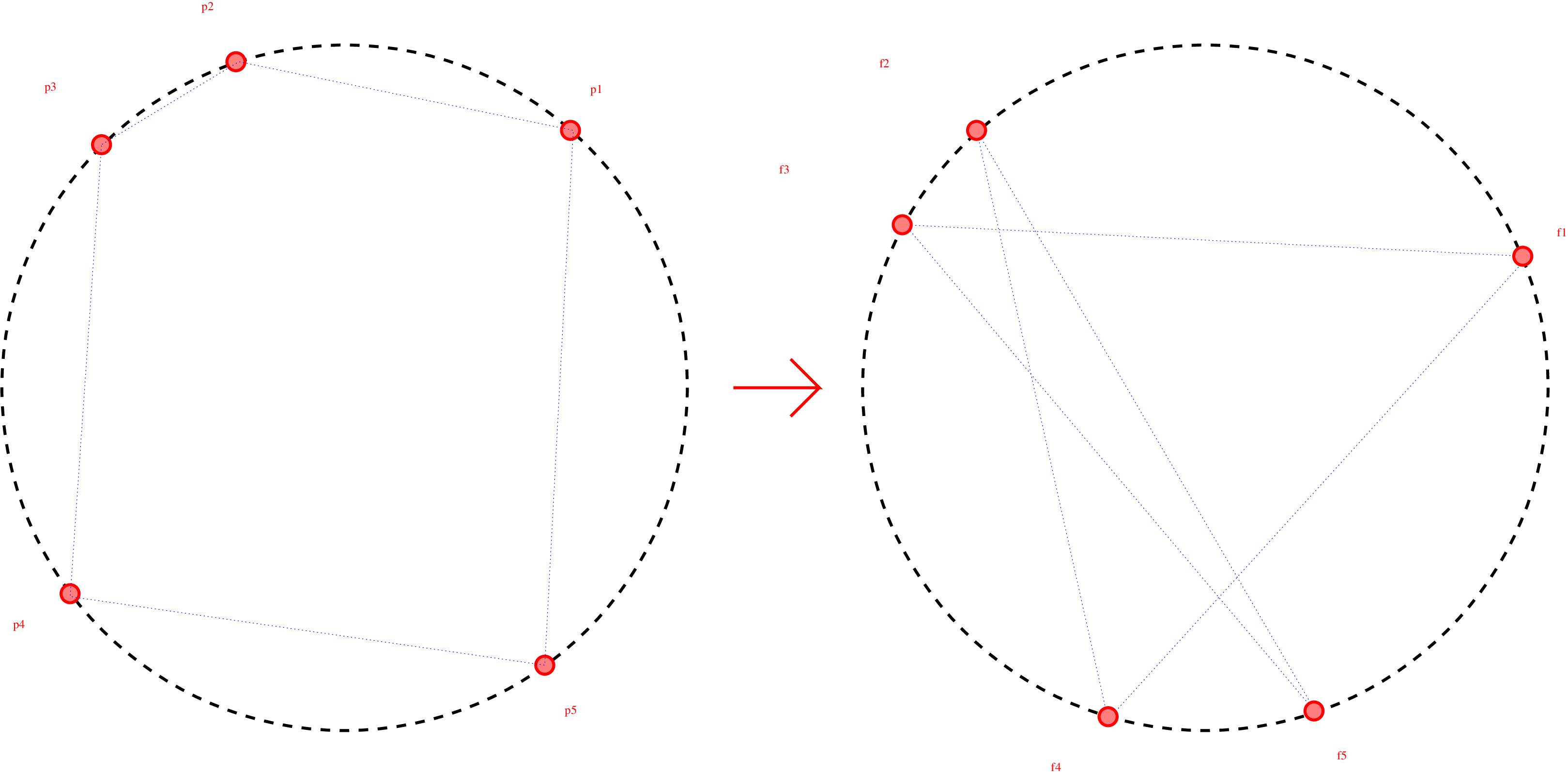}$$
\centerline{\hskip 0.6cm cyclic order \hskip 2.3cm non-consecutive order}
\caption{Five points in standard cyclic order in $S^2$ mapped to points in non-consecutive order in $S^4$.}
}\label{pentfig}
\end{figure}

By non-consecutive order we mean that if $C \subset S^4$ is the round circle such that $f_*(p) \in C_5(C)$ and if
$I \subset C$ is an embedded interval such that $\partial I = \{ f(p_i), f(p_{i+1}) \}$ then the interior 
of $I$ must contain a third point of $f_*(p)$.  For this purpose our indices $i$ are taken modulo $5$, $i \in \Zed_5$. 
See Figure \ref{pentfig}.  Similarly, points $p$ being in the cyclic order of the circle they lie in means that for 
any $i$ one can choose an embedded interval $I \subset C$ such that $I$ contains only the two points
$\{p_i, p_{i+1}\}$ of $p$.

\begin{defn}\label{condef} Given a space $M$ the configuration-space of distinct $k$-tuples of points in $M$ is defined as 
$$C_k(M) = \{ (p_1, \cdots, p_k) \in M^k : p_i \neq p_j \ \forall i \neq j \}.$$ 

We think of $S^n$ as the unit sphere in $\Real^{n+1}$.  A round circle in $S^n$ is the intersection of 
$S^n$ with a $2$-dimensional affine-linear subspace of $\Real^{n+1}$. 
\end{defn}

As we will see, generically the set $\mathcal C_f$ is a canonically-oriented $2$-dimensional compact manifold.  If one composes the
inclusion $\mathcal C_f \to C_5(S^2)$ with any of the forgetful maps $O_i : C_5(S^2) \to S^2$ where 
$O_i(p_1, \cdots, p_5) = p_i$ this gives a map between compact, oriented $2$-dimensional manifolds
$\mathcal C_f \to S^2$.  One could ask, what is the degree of this map?  We will see it is always zero due
to a symmetry issue. Notice that if $p \in \mathcal C_f$ then the reversal $\overline{p}$ is also an element of $\mathcal C_f$ where
$\overline{p} = (p_4, p_3, p_2, p_1, p_5)$ (the reversal that fixes $p_5$).   Think of the reversal operation
as an involution of $C_5 M$, for any $M$.  Given that reversal is fixed-point free, there is an induced map of compact manifolds
$\mathcal C_f / \Zed_2 \to S^2$ where we are now forced to compose the inclusion $\mathcal C_f / \Zed_2 \to C_5(S^2)/\Zed_2$ 
with the forgetful map $O_5$.  We will see that
$\mathcal C_f / \Zed_2$ is naturally only a compact $2$-manifold, i.e. it does not inherit a canonical orientation, thus we can
only take the mod-$2$ degree of this map.   This invariant is the topic of our paper. 

\begin{thm}\label{mainthm}
Given a smooth $2$-knot $f : S^2 \to S^4$ the mod-$2$ degree of the map $\mathcal C_f / \Zed_2 \to S^2$ 
as an element of $\Zed_2$ is an isotopy-invariant of $f$.   As a map
of the form $\mu : \pi_0 \Emb(S^2, S^4) \to \Zed_2$ it is additive with respect to the connect-sum monoid structure on the domain.
Moreover there are knots with $\mu \neq 0$.
\end{thm}

Thus the invariant $\mu$ can be thought of as a coarse measure of the extent to which smooth embeddings `shuffle' five points 
on a round circle. 

Stereographic projection at a point $p \in S^n$ can be thought of as a map $S^n \to (T_p S^n) \cup \{\infty\}$.  From this
perspective it is a conformal diffeomorphism. Such conformal diffeomorphisms are known to send round circles in $S^n$ to the either 
round circles or straight lines in $T_p S^n$, depending on whether or not the circle runs through the stereographic projection point $p$. 
Similarly, if one
stereographically projects a $2$-knot $f$ at some point in its image, it converts the $2$-knot into what is commonly called a {\bf long knot}. 
This leads to a homotopy-equivalence \cite{BudSurv} $\Emb(S^2, S^4) \simeq SO_5 \times_{SO_2} \Emb(D^2, D^4)$, where 
$\Emb(D^2, D^4)$ denotes the space of smooth embeddings $g : \Real^2 \to \Real^4$ such that $g(D^2) \subset D^4$ and
$g(p) = (p,0)$ for all $p \in \Real^2 \setminus D^2$.  Notice that up to an isometry of $S^4$, $g$ is precisely the sterographic 
projection of a knot $f : S^2 \to S^4$ which is a standard linear embedding on a hemisphere. 

Take the perspective of computing $\mu(f)$ by counting points in the pre-image of a regular value $p \in S^2$ for
the map $\mathcal C_f / \Zed_2 \to S^2$. Our homotopy-equivalence
above gives us a natural isomorphism $\pi_0 \Emb(D^2, D^4) \to \pi_0 \Emb(S^2, S^4)$, thus we can define $\mu : \pi_0 \Emb(D^2, D^4) \to \Zed_2$
as $\mu$ of its stereographic projection in $\Emb(S^2, S^4)$.  Thus if $g \in \Emb(D^2, D^4)$ then $\mu(g)$ is a count of the 
linearly-ordered quadrisecants in $D^2$ that are mapped by $f$ to an alternating quadrisecant in $D^4$.  Specifically, $\mu(g)$ is the count of 
points (up to reversal) $p \in C_4 D^2$ that sit on an oriented straight line $L \subset \Real^2$ with $p_1 < p_2 < p_3 < p_4$ in the linear 
order of $L$ such that $g(p)$ also sits on an oriented straight line 
$L' \subset \Real^4$ with $g(p_3) < g(p_1) < g(p_4) < g(p_2)$ in the linear order of $L'$. This is our approach to proving Theorem 
\ref{mainthm}.  By applying some standard techniques in the homology of configuration spaces we
turn these formulas into something analogous to Polyak-Viro formulas, which allows us to compute the invariant Whitney diagrams
or Yoshikawa diagrams of $2$-knots. 

Further, we extend these arguments to show there is an invariant defined for families of knots

$$\mu : \pi_{2(n-j-2)} \Emb(S^j, S^n) \to 
\begin{cases} \Zed & \text{if } j=1 \text{ or both } n \text{ and } j \text{ odd}\\ 
              \Zed_2 & \text{otherwise.}\end{cases}$$

This invariant is similarly defined by considering $5$-tuples on a round circle in $S^j$ that are mapped
to $5$ points that sit on a round circle in $S^n$ (through the entire family), with cyclic ordering 
being sent to non-consecutive ordering.  One mods-out by the reversal involution
(this step is replaced when $j=1$ by restricting to the counter-clockwise cylically-ordered subspace of $C_5 S^1$) and takes the degree or 
mod-$2$ degree of the forgetful map, as appropriate.  
This homomorphism is defined whenever both inequalities $n \geq j+2$ and $j \geq 1$ hold.   

The invariant $\mu$ has been well-studied when $j=1$.  In the paper \cite{BCSS} the authors observed
for $(n,j)=(3,1)$ that $\mu$ equals the type-$2$ invariant of knots. The invariant was written-up in the quadrisecant form
in \cite{BCSS} while the M.Sc thesis of Flowers \cite{Flowers} put it in the language of circular pentagrams. 
See the \href{https://sean564.github.io/top/}{demonstration} by Sean Lee. 
In the follow-up paper \cite{BudSurv} it was observed for the $(n,j)$ case with $n \geq 4$ and $j=1$ that 
the homomorphism $\mu : \pi_{2n-6} \Emb(D^1, D^n) \to \Zed$ is an isomorphism of groups. Moreover it was known at the time that the 
lowest-dimensional non-trivial homotopy
group of the space $\Emb(D^1, D^n)$ was in dimension $2n-6$.  The first non-trivial homotopy group of $\Emb(D^j, D^n)$ 
is known to occur in dimension $2n-3j-3$ when $2n-3j-3 \geq 0$, and otherwise it is typically in dimension zero \cite{BudSurv}. 
That said there remains some important open cases such as the question of the triviality of $\pi_0 \Emb(D^4, D^4) = \pi_0 \Diff(D^4)$. 

This paper was inspired by the sequence of papers \cite{BG}, \cite{BG2}, \cite{BG3}, where analogous
invariants were defined out of groups such as $\pi_{2n-6} \Emb(I, S^1 \times D^{n-1})$.   The main result of this paper
was stumbled-upon during a visit to Princeton, while preparing a presentation.  The author would like to thank 
David Gabai for hosting, and Scott Carter for his early comments, as well as Danny Ruberman, Victor Turchin and Tadayuki Watanabe
for their comments on an early draft. 

\section{The invariant}\label{invtsec}

We begin with a detailed description of the $\mu$ invariant in the form

$$\mu : \pi_{2(n-j-2)} \Emb(D^j, D^n) \to 
\begin{cases} \Zed & \text{if } j=1 \text{ or both } n \text{ and } j \text{ odd}\\ 
              \Zed_2 & \text{otherwise.}\end{cases}$$

We will use the language of intersections of maps with submanifolds (transversal intersections).  
Let $\mathcal L^s \subset C_4(D^j)$ be the subspace defined by the condition on $p \in C_4(D^j)$ 
that there exists an oriented straight line $L \subset \Real^j$ such that 
$p \in C_4(L)$ and $p_1 < p_2 < p_3 < p_4$ in the linear ordering on $L$ induced from its orientation.  Similarly, 
we let $\mathcal L^a \subset C_4(D^n)$ be the subset of points $p \in C_4(D^n)$ such that there exists an oriented line $L$
with $p_3 < p_1 < p_4 < p_2$ with respect to the linear order on $L$ induced from its orientation.  Given 
$f : S^{2(n-j-2)} \to \Emb(S^j, S^n)$ let $f_* : S^{2(n-j-2)} \times C_4 D^j \to C_4 D^n$ be the induced map, i.e.
$f_*(v, p_1, \cdots, p_4) = (f(v)(p_1), f(v)(p_2), f(v)(p_3), f(v)(p_4))$.  
Let $i : \mathcal L^s \to C_4(D^j)$ be the inclusion, and $I$ the identity map on $S^{2(n-j-2)}$.  Provided 
$f_* \circ (I \times i) : S^{2(n-j-2)} \times \mathcal L^s \to C_4(D^n)$ is transverse to $\mathcal L^a$, 
$\mu(f)$ will be defined in terms of this intersection.    As in the papers \cite{BCSS} \cite{Flowers} one can 
apply a small perturbation to the map $f : S^{2(n-j-2)} \to \Emb(S^j, S^n)$ to ensure transversality of the
family $f_*$.  But for the purpose of the definition here would could simply perturb $f_* \circ (I \times i)$ to be transverse.  

Let $R : C_4(D^k) \to C_4(D^k)$ be the reversal involution $R(p_1,p_2,p_3,p_4) = (p_4,p_3,p_2,p_1)$.  The map $R$ 
restricts to an involution of $\mathcal L^s$ and $\mathcal L^a$, respectively, making $f_* \circ (I \times i)$ an
equivariant map. 

\begin{lem}
The involution $R$ of $C_4(D^k)$ is orientation-preserving for all $k$, interpreting $C_4(D^k)$ as an open subset of $(D^k)^4$
with its standard product orientation. When restricted to $\mathcal L^s$ (in $C_4(D^j)$) it multiplies the orientation by $(-1)^{j+1}$.  Similarly, 
it multiplies the orientation of $\mathcal L^a$ (in $C_4(D^n)$) by $(-1)^{n+1}$. 
\end{lem}

\begin{defn}
Our invariant $\mu(f)$ is defined as the signed intersection number of 
$f : S^{2(n-j-2)} \times \mathcal L^s / \Zed_2 \to C_4(D^n)/\Zed_2$ with the subspace $\mathcal L^a/\Zed_2$, if both 
manifolds are oriented and $j > 1$.  If either manifold fails to be oriented and $j>1$ we use the mod-$2$ intersection number.
In the special case of $j=1$, $\mathcal L^s_+$ has two path-components, we let $\mathcal L^s_+$ denote the 
component where the points are in increasing order (in the linear ordering of $\Real$), and define $\mu(f)$ as the signed intersection number of 
$f : S^{2(n-j-2)} \times \mathcal L^s_+ \to C_4(D^n)$ with $\mathcal L^a$. 
\end{defn}

\begin{prop}
$$\mu : \pi_{2(n-j-2)} \Emb(D^j, D^n) \to 
\begin{cases} \Zed & \text{if } j=1 \text{ or both } n \text{ and } j \text{ odd}\\ 
              \Zed_2 & \text{otherwise}\end{cases}$$
is well-defined. 
\begin{proof}
While the manifolds $C_4(D^j)$, $C_4(D^n)$, $\mathcal L^a$ and $\mathcal L^s$ are non-compact, the transverse intersection
of $f_* : S^{2(n-j-2)} \times \mathcal L^s \to C_4(D^n)$ with $\mathcal L^a$ is compact, and given a homotopy 
$H : I \times S^{2(n-j-2)} \times \mathcal L^s \to C_4(D^n)$ the (transverse) intersection of $H$ with $\mathcal L^a$ is 
also compact.  

There is an essentially analytic argument for this. Given a smooth embedding $g \in \Emb(D^j, D^n)$ there is a lower bound 
$\epsilon > 0$ on how close points can be in any quadrisecant. An application of the triangle inequality 
shows that $g$ satisfies a reverse-Lipschitz inequality 
$$(m-KR) \cdot |x-y| \leq |g(x)-g(y)|$$
provided $|x-y| < R$.  In the inequality, $K$ is a Lipschitz constant for $g'$, i.e. $||g'_x-g'_y|| \leq K |x-y|$ for all $x,y \in D^j$
and $m = \min \{ |g'_p(v)| : p \in D^j, v \in S^{j-1} \}$. Thus $\epsilon = \frac{m}{2K}$ gives the reverse-Lipschitz inequality, and
thus if points are closer than $\epsilon$ their linear ordering (on any line) is preserved.  Thus the constants $m, K$ and $R$ can be chosen
continuously for a $C^2$-family $f$.
\end{proof}
\end{prop}

While the geometry of this invariant is appealing -- it literally is a measure of how embeddings `shuffle' quadruples of points 
along straight lines -- it leaves us with the problem of how such an invariant can be practically computed.

\subsection{A double-point formulation of the invariant, the $P_\epsilon$ deformation.}

Given that the invariant is an intersection number, it is essentially homological in nature -- in the homology of
configuration spaces.  This gives us considerable flexibility in the computation of the invariant.  The computation here is
largely inspired by \cite{BG} and \cite{BG2},  specifically Lemma 3.4 in \cite{BG2}. That lemma was in turn inspired by a 
(yet unpublished) argument of Misha Polyak's describing a clean relation between this invariant in the $(n,j)=(3,1)$ case \cite{BCSS} 
and the Polyak-Viro perspective on the type-$2$ invariant \cite{PV98}.  These arguments could also be considered as flowing from 
the perspective of the Gravity Filtration popularized by Fred Cohen (Reference \cite{Tam} is a good example of the idea). 
One other similar thread of ideas follows from the integral expressions for Vassiliev invariants involving volume forms on 
spheres -- in principle our formulas could be viewed as the limiting expressions as one perturbs the volume forms to delta functions.

For $\epsilon \in \Real$ consider the diffeomorphism $P_\epsilon$ of $\Real^n$ given by 
$$P_\epsilon(x_1, x_2, \cdots, x_n) = \left(x_1, x_2, \cdots, x_{n-1}, x_n + \epsilon \sum_{i=1}^{n-1} x_i^2\right).$$
This diffeomorphism has the feature that it converts the $x_n=c$ hyperplanes into paraboloids, when $\epsilon \neq 0$. Similarly 
it turns lines into parabolas, with the exception that it acts by translation on the lines parallel to the $x_n$-axis.
Moreover this is a group action of $\Real$ on $\Diff(\Real^n)$.  Lines parallel to the $x_n$-axis we will call {\it vertical} and
two or more points on a common vertical line we will similarly call vertical pairs, triples, quadruples, etc. 

\begin{figure}[H]
{
\psfrag{zeta}[tl][tl][0.7][0]{$x_n$}
\psfrag{2}[tl][tl][0.7][0]{$2$}
\psfrag{3}[tl][tl][0.7][0]{$3$}
\psfrag{km3}[tl][tl][0.7][0]{$k-3$}
\psfrag{km2}[tl][tl][0.7][0]{$k-2$}
\psfrag{km1}[tl][tl][0.7][0]{$k-1$}
$$\includegraphics[width=14cm]{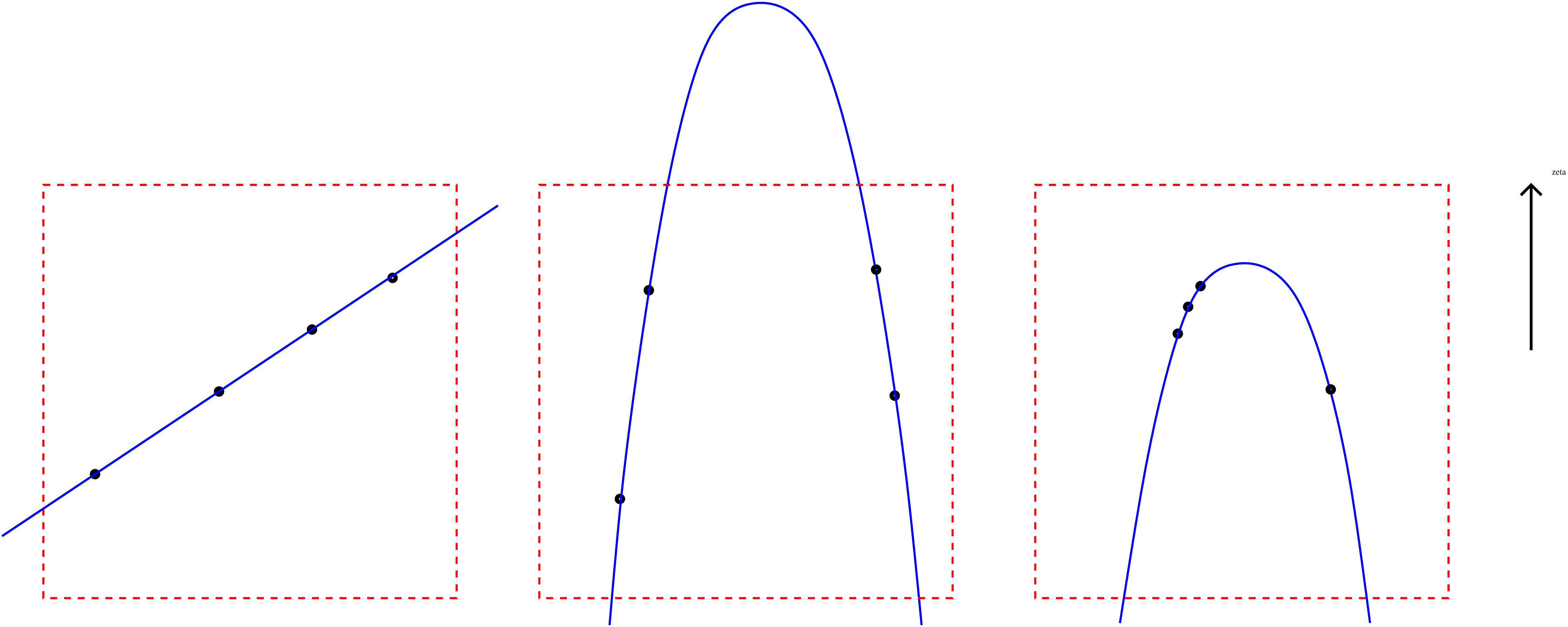}$$
\caption{Parabolic quadruples, $\epsilon=0$ left.  Large $\epsilon$ middle and right.}\label{parabdiag}
}
\end{figure}

The motivation for introducing $P_\epsilon$ is that rather than computing the invariant $\mu$ using the original families
of quadruples, $\mathcal L^a$ and $\mathcal L^s$, we use their images $P_\epsilon^*(\mathcal L^a)$ and
$P_\epsilon^*(\mathcal L^s)$, as depicted in Figure \ref{parabdiag}, i.e. in our family of maps
$S^{2(n-j-2)} \times D^j \to D^n$ we would be counting $4$-tuples of points on the appropriate parabola in
$D^j$ being mapped to points on the appropriate parabola in $D^n$.  Our interest comes from observing how these computations
trend as $\epsilon \to \infty$. 

For the remainder of this section we restrict to the $(n,j) = (4,2)$ case.  When we project a $2$-knot of the form 
$f : S^2 \to \Real^4$ (or $f \in \Emb(D^2, D^4)$) into a $3$-dimensional vector subspace of 
$\Real^4$ (for $f \in \Emb(D^2, D^4)$ this
subspace should contain the long axis), this map can generically be assumed to be locally 
an immersion at all but finitely many points, and those finite points are called `cross-caps' or `Whitney umbrellas'
\cite{Whit}. There will be a $1$-manifold of double points, a $0$-manifold of triple points, and no quadruple points. 
The cross-caps are not isolated from the double-point curves, as double-point curves 
can terminate at cross-caps.   These observations will help us compute $\mu$ of a $2$-knot. It turns out cross-caps can be removed
via an isotopy \cite{CG}, although we will not use this.  

Given a parabola of $P^*_\epsilon(\mathcal L^a)$ intersecting a $2$-knot, we can assume the points of intersection do not
include the maximum ($x_4$-coordinate) of the parabola, as such parabolas (generically) approximate vertical quadruples on
the knot, which generically do not occur.  Thus the four points of the parabola intersect the knot are partitioned in two
groups, determined by which side of the maximum they are on.  There are three possibilities, $4+0$, $3+1$ and $2+2$.  The
$4+0$ case can not occur as it corresponds to a vertical quadruple, which generically does not exist.  Thus we have
only the latter two possibilities.  In the limit they come from intersection with the submanifolds of $C_4(D^4)$ below,
with the $3+1$ case occuring in two variations.
Note that when a $2$-knot $S^2 \to \Real^4$
is in Whitney's general-position with respect to a projection map $\Real^4 \to \Real^3$, the induced map
$C_4 S^2 \to C_4 \Real^4$ is transverse to the double-point and triple-point submanifolds.  In particular the double and
triple-points can be thought of as a stratified subset of the domain $S^2$, given by a collection of isolated points 
corresponding to the Whitney umbrellas and triple points, and a collection of $1$-manifolds corresponding to the double-point
curves. 

Given that the family $P^*_\epsilon(\mathcal L^a)$ is not an isotopy on the interval $[0,\infty]$ (it is an isotopy on $[0,\infty)$), one needs to further check that in a neighbourhood of $\epsilon=\infty$ that double-point pairs and triple points have essentially
unique parabolic quadruples, for $\epsilon$ finite.  This rapidly follows from the geometry of parabolas when applied to
maps in Whitney's general position, much like in the case of classical knots and the corresponding Polyak-Viro formulas.

\begin{figure}[H]
{
\psfrag{p1}[tl][tl][0.7][0]{$p_1$}
\psfrag{p2}[tl][tl][0.7][0]{$p_2$}
\psfrag{p3}[tl][tl][0.7][0]{$p_3$}
\psfrag{p4}[tl][tl][0.7][0]{$p_4$}
\psfrag{v1}[tl][tl][0.7][0]{\textcolor{green}{$v$}}
\psfrag{v2}[tl][tl][0.7][0]{\textcolor{green}{$w$}}
\psfrag{v3}[tl][tl][0.7][0]{\textcolor{green}{$\alpha'(v)$}}
\psfrag{v4}[tl][tl][0.7][0]{\textcolor{green}{$\beta'(w)$}}
\psfrag{L}[tl][tl][1][0]{\textcolor{red}{$L$}}
$$\includegraphics[width=14cm]{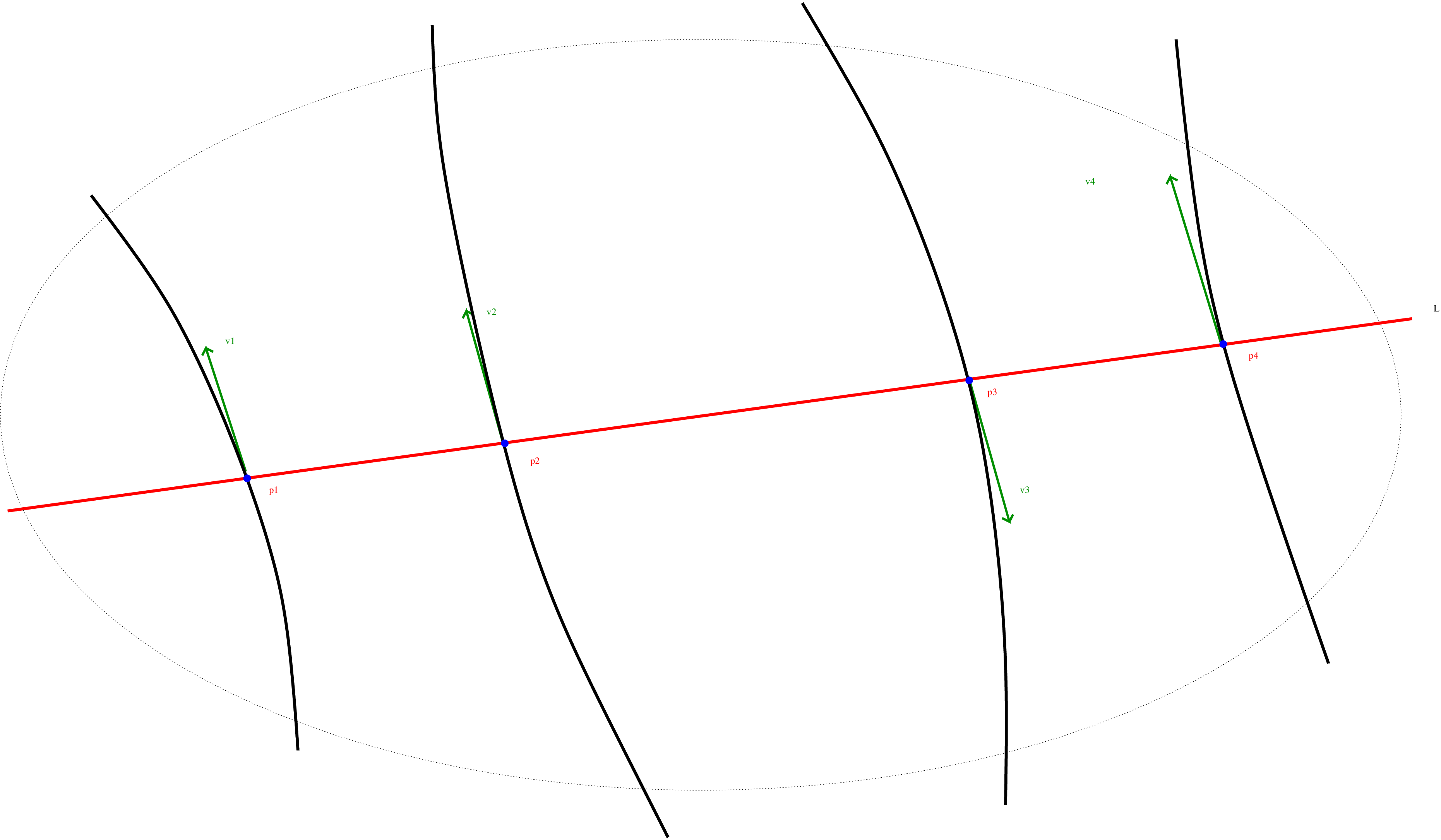}$$
\caption{Transversality for quadrisecants intersecting double-point manifolds}\label{tqd_trans}
}
\end{figure}

Consider the issue of quadrisecants (thought of as a submanifold of $C_4 D^2$ transversely intersecting the double and
triple point submanifolds of $C_4 D^2$. Triple points and Whitney umbrellas can be avoided, as a small diffeomorphism
of the domain $D^2$ allows us to make a small perturbation of a point (and generically three points are not on a line).
Thus our transversality condition can be stated purely in terms of the collinear quadruple manifold where the points
$p_1 < p_2 < p_3 < p_4$ on the line $L$ in that linear order intersect the double-point manifold, i.e. 
$f(p_3)$ is under $f(p_1)$ and $f(p_2)$ is under $f(p_4)$. Let $\alpha$ be the diffeomorphism of a neighbourhood of
the double-point curve at $p_1$ that associates to it double-point in the neighbourhood of $p_3$ corresponding to the
overcrossing.  Similarly, let $\beta$ be the diffeomorphism of the neighbourhood of $p_4$ in the double-point curve
such that $f(p_2)$ is under $f(p_4)$.  If we parametrize the quadrisecant manifold via the map
$C_2(\Real^2) \times (0,\infty)^2 \to C_4(\Real^2)$ by mapping
$(p,q, t, s) \longmapsto (p,q, q+t(q-p), q+(t+s)(q-p)) = (p_1,p_2,p_3,p_4)$ then the transversality of our intersection
corresponds to the non-degeneracy of the $8 \times 8$-matrix

$$\begin{pmatrix} 
v & 0 & I & 0 & 0 & 0 \\
0 & w & 0 & I & 0 & 0 \\
\alpha'(v) & 0 & (1-t)I & tI & q-p & 0 \\
0 & \beta'(w) & (1-t-s)I & (t+s)I & q-p & q-p\end{pmatrix}.$$

As described in
Figure \ref{tqd_trans}, the vectors $v$ and $w$ are tangent vectors to the double-point curves at $p_1=p$
and $p_2=q$.
Most of the entries in the above matrix are $2 \times 1$ vectors, such as $v$.  But the sub-matrix $I$ is the $2 \times 2$
identity matrix. With a little algebra one can show this is equivalent to the non-degeneracy of the matrix
$$\begin{pmatrix}
\alpha'(v)  + (t-1)v & -tw & q-p & 0 \\
(t+s-1)v & \beta'(w) -(t+s)w & 0 & q-p \end{pmatrix}.$$
Roughly speaking, this expression being equal to zero is equivalent to there being an (infinitesimal)
$1$-parameter family of quadrisecants through the double-point curves, or the tangent vectors to
our curves being parallel to the quadrisecant direction.

Let $\mathcal{C}_{\ctwo{1}{3}\ctwo{4}{2}}$ be the subspace of $C_4(D^4)$ such that $p_1$ is over $p_3$ and $p_4$ is over $p_2$ with respect
to the coordinate projection $\Real^4 \to \Real^3$ given by $(x_1,x_2,x_3,x_4) \longmapsto (x_1,x_2,x_3)$. In 
summary, we have the intermediate result.

\begin{prop}\label{halfprop}
Given a long knot $f : D^2 \to D^4$, provided the induced map $f_* : C_4 D^2 \to C_4 D^4$ is such that
$\mathcal L^s \subset C_4 D^2$ is transverse to $f_*^{-1}(\mathcal{C}_{\ctwo{1}{3}\ctwo{4}{2}})$ then 
$\mu(f)$ is the mod two count of the points in 
$\left( \mathcal L^s \cap f_*^{-1}(\mathcal{C}_{\ctwo{1}{3}\ctwo{4}{2}}) \right)/\Zed_2$
where the involution $\Zed_2$ is the reversal $(p_1,p_2,p_3,p_4) \longmapsto (p_4,p_3,p_2,p_1)$ action
on $C_4(D^2)$.
\end{prop}

In Subsection \ref{ss2} we will manipulate the formula of Proposition \ref{halfprop} to make it 
somewhat easier to work with, in practical computations.  One quick observation that can be derived from
Proposition \ref{halfprop} appears below. In Subsection \ref{ss2} we will have a formula for $\mu$
that is perhaps the direct analogue to a Polyak-Viro formula.  

\begin{cor}\label{artinspuncor}
If $f : D^2 \to D^4$ is Artin-spun from a classical knot $D^1 \to D^3$ then $\mu(f)=0$.
\begin{proof}
For an Artin-spun knot, the double-point diagram consists of a collection of nested 
concentric circles in $D^2$, moreover on any ray out of the origin the double-point diagram
restricts to the chord diagram for the original knot $D^1 \to D^3$.  Thus the 
quadrisecants of $\mathcal L^s \cap f_*^{-1}(\mathcal{C}_{\ctwo{1}{3}\ctwo{4}{2}})$ will precisely be the quadrisecants
for the original knot, but by the circular symmetry it will be $SO_2$-invariant.  In particular
this is a non-transverse intersection.  For each quadrisecant, we can eliminate its $SO_2$-invariant
double-point curve from the diagram by pre-composing our embedding $f : D^2 \to D^4$ by a small diffeomorphism
supported in a neighbourhood of the double-point curve -- apply a slight $SO_2$-motion along the curve, damping
it out in the neighbourhood. One by one this eliminates all the double-point curves.
\end{proof}
\end{cor}

\subsection{Applying $P_\epsilon$ to the quadrisecants in the domain.}\label{ss2}

We now consider applying the same family of diffeomorphisms $P_\epsilon$ to the quadrisecants in the domain $D^2$, 
i.e. we modify the $\mathcal L^s$ manifold in the count
$$\left( \mathcal L^s \cap f_*^{-1}(\mathcal{C}_{\ctwo{1}{3}\ctwo{4}{2}}) \right)/\Zed_2$$
by considering
$$\left( P^*_\epsilon(\mathcal L^s) \cap f_*^{-1}(\mathcal{C}_{\ctwo{1}{3}\ctwo{4}{2}})) \right) / \Zed_2$$
with $\epsilon$ large.   The naive limit would be the intersection

$$\left( \left( \mathcal{C}_{\cfour{4}{3}{2}{1}} \cup 
\mathcal{C}_{\cthree{3}{2}{1}\cone{4}} \cup 
\mathcal{C}_{\ctwo{2}{1}\ctwo{3}{4}} \cup 
\mathcal{C}_{\cone{1}\cthree{2}{3}{4}} \cup
\mathcal{C}_{\cfour{1}{2}{3}{4}}
\right) \cap 
f_*^{-1}(\mathcal{C}_{\ctwo{1}{3}\ctwo{4}{2}}) \right) / \Zed_2. \hskip 1cm (*)$$

Specifically, we ask if this generically is the same count as the intersection
$\left( \mathcal L^s \cap f_*^{-1}(\mathcal{C}_{\ctwo{1}{3}\ctwo{4}{2}}) \right)/\Zed_2$.

The set $\mathcal{C}_{\cthree{3}{2}{1}\cone{4}}$ denotes the subspace of $C_4(D^2)$ where
$p_3$ is over $p_2$ which is in turn over $p_1$, in that order, thus it is a co-dimension $2$
submanifold of $C_4(D^2)$.  Similarly, $\mathcal{C}_{\cfour{4}{3}{2}{1}}$ denotes the subspace of a
vertical quadrisecant, in the given order, thus this has co-dimension $3$.  In particular
the intersection $\left( \mathcal{C}_{\cfour{4}{3}{2}{1}} \cup 
\mathcal{C}_{\cfour{1}{2}{3}{4}}
\right) \cap 
f_*^{-1}(\mathcal{C}_{\ctwo{1}{3}\ctwo{4}{2}})$ generically is empty. So (*) generically is given by the intersection

$$\left( \left( \mathcal{C}_{\cthree{3}{2}{1}\cone{4}} \cup 
\mathcal{C}_{\ctwo{2}{1}\ctwo{3}{4}} \cup 
\mathcal{C}_{\cone{1}\cthree{2}{3}{4}} 
\right) \cap 
f_*^{-1}(\mathcal{C}_{\ctwo{1}{3}\ctwo{4}{2}}) \right) / \Zed_2. \hskip 1cm (*)$$

In this section we will compare this count to the count 
$\left( P^*_\epsilon(\mathcal L^s) \cap f_*^{-1}(\mathcal{C}_{\ctwo{1}{3}\ctwo{4}{2}})) \right) / \Zed_2$ 
for $\epsilon$ large, given that the transition is no longer from an isotopy of manifolds, there is the
issue that parabolic quadruples could collapse (multiple-to-one) that we need to address. 

\begin{figure}[H]
{
\psfrag{a}[tl][tl][0.7][0]{(a)}
\psfrag{b}[tl][tl][0.7][0]{(b)}
\psfrag{p1}[tl][tl][0.7][0]{$p_1$}
\psfrag{p2}[tl][tl][0.7][0]{$p_2$}
\psfrag{p3}[tl][tl][0.7][0]{$p_3$}
\psfrag{p4}[tl][tl][0.7][0]{$p_4$}
$$\includegraphics[width=8cm]{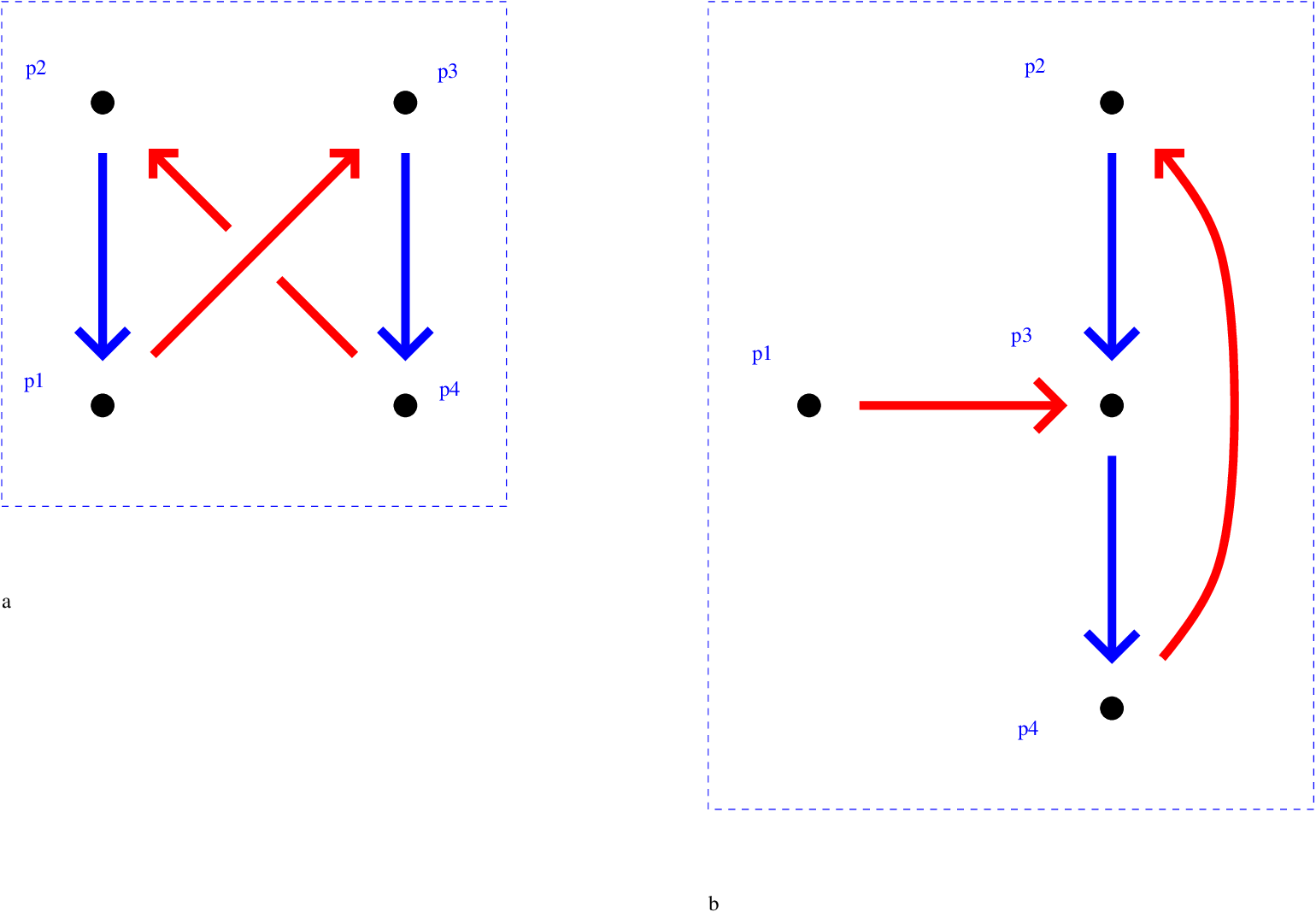}$$
\caption{The invariant $\mu$ as a count of vertical tuples.}\label{cd0vt}
}
\end{figure}

Figure \ref{cd0vt} gives a diagrammatic description of elements in the intersection of type $(*)$.  
There is one diagram for each $R$-orbit ($R$ the reversal involution). 
A blue arrow $p \textcolor{blue}{\to} q$ means the point $p$ is over the point $q$ in the 
domain of $f : D^2 \to D^4$.  A red arrow $p \textcolor{red}{\to} q$ means $f(p)$ is over the point $f(q)$ in the codomain of $f$.  
Thus Figure \ref{cd0vt} (a) indicates that there is a `$4$-cycle of overcrossings' in the sense that:  
\begin{itemize}
\item $p_2$ is over $p_1$ in $D^2$, 
\item $f(p_1)$ is over $f(p_3)$ in $D^4$, 
\item $p_3$ is over $p_4$ in $D^2$, and finally 
\item $f(p_4)$ is over $f(p_2)$ in $D^4$.  
\end{itemize}

Figure \ref{cd0vt} (a) depicts elements of 
$$\mathcal{C}_{\ctwo{2}{1}\ctwo{3}{4}} \cap f_*^{-1}(\mathcal{C}_{\ctwo{1}{3}\ctwo{4}{2}})$$

Similarly, Figure \ref{cd0vt} (b) $$\left(  %\mathcal{C}_{\cthree{3}{2}{1}\cone{4}} \cup 
\mathcal{C}_{\cone{1}\cthree{2}{3}{4}} \cap 
f_*^{-1}(\mathcal{C}_{\ctwo{1}{3}\ctwo{4}{2}}) \right) / \Zed_2$$

is the pre-image of (1) intersect a $3+1$ decomposition in the domain.
This could be described as a `$2$-cycle of overcrossings' in the sense that $p_2$ is over $p_4$ in $D^2$
while $f(p_4)$ is over $f(p_2)$ in $D^4$, with the exception of the requirement of the intermediate point $p_3$, 
between $p_2$ and $p_4$ in $D^2$, which is itself an undercrossing in the sense that $f(p_3)$ is under $f(p_1)$ in $D^4$. 

Figure \ref{cd0vt} (b) describes elements of intersection type 
$\left(  %\mathcal{C}_{\cthree{3}{2}{1}\cone{4}} \cup 
\mathcal{C}_{\cone{1}\cthree{2}{3}{4}} \cap 
f_*^{-1}(\mathcal{C}_{\ctwo{1}{3}\ctwo{4}{2}}) \right) / \Zed_2$
with the reverse being the intersection of type $\left(  \mathcal{C}_{\cthree{3}{2}{1}\cone{4}}  
%\mathcal{C}_{\cone{1}\cthree{2}{3}{4}} 
\cap  f_*^{-1}(\mathcal{C}_{\ctwo{1}{3}\ctwo{4}{2}}) \right) / \Zed_2.$

\begin{figure}[H]
{
\psfrag{p0}[tl][tl][0.7][0]{$p_1$}
\psfrag{p1}[tl][tl][0.7][0]{$p_2$}
\psfrag{p2}[tl][tl][0.7][0]{$p_3$}
\psfrag{p3}[tl][tl][0.7][0]{$p_4$}
\psfrag{I0}[tl][tl][0.7][0]{$I_1$}
\psfrag{I1}[tl][tl][0.7][0]{$I_2$}
\psfrag{I2}[tl][tl][0.7][0]{$I_3$}
\psfrag{I3}[tl][tl][0.7][0]{$I_4$}
$$\includegraphics[width=8cm]{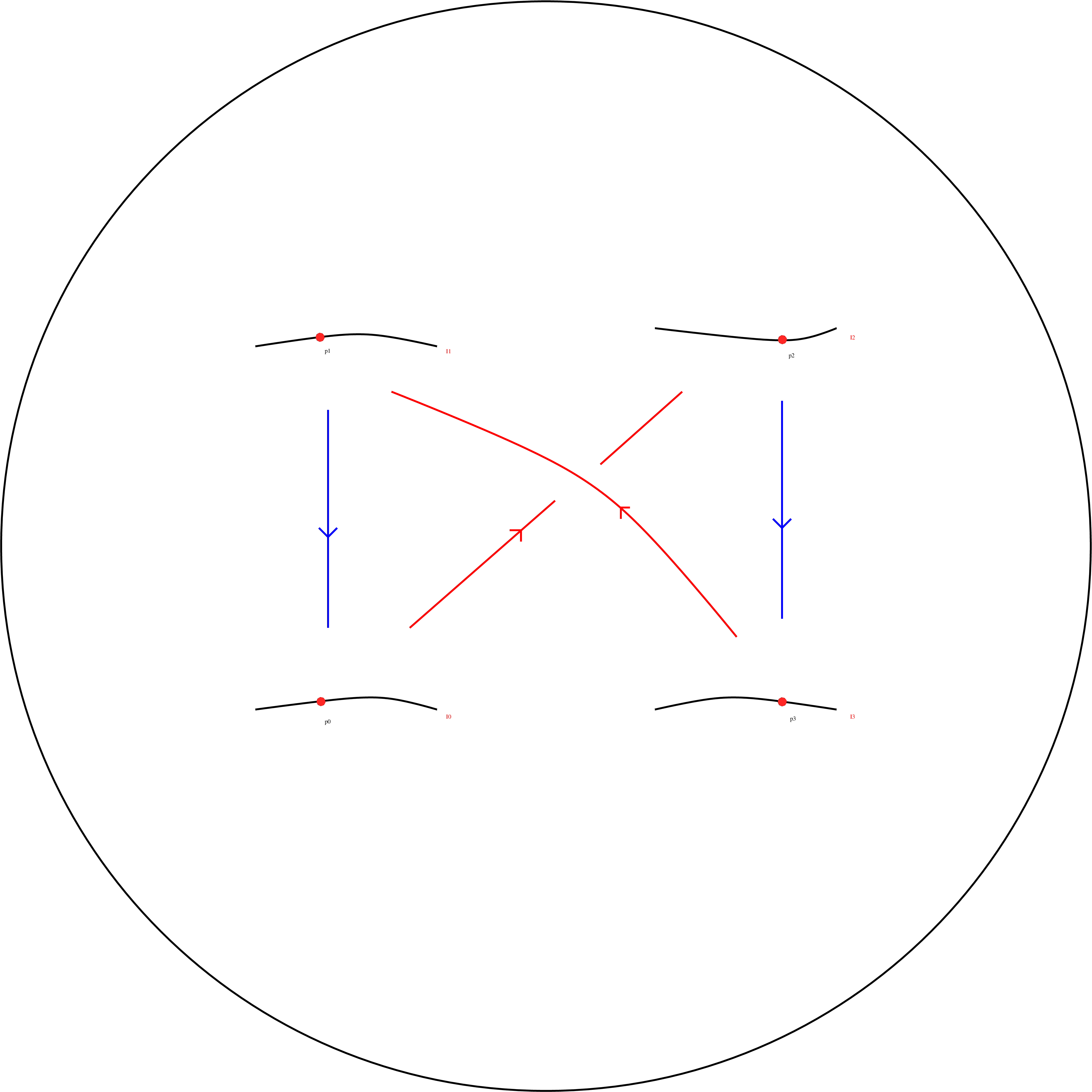}$$
\caption{Transversality for $\mathcal{C}_{\ctwo{2}{1}\ctwo{3}{4}} \cap f_*^{-1}(\mathcal{C}_{\ctwo{1}{3}\ctwo{4}{2}})$}\label{trans1fig}
}
\end{figure}

To begin, consider the transversality condition for the intersections of Figure \ref{cd0vt} (a)
$$\mathcal{C}_{\ctwo{2}{1}\ctwo{3}{4}} \cap f_*^{-1}(\mathcal{C}_{\ctwo{1}{3}\ctwo{4}{2}}).$$

depicted in Figure \ref{trans1fig}. If we let $\alpha : I_1 \to I_3$ be the double-point 
diffeomorphism from a neighbourhood of $p_1$ to a neighbourhood of $p_3$, and $\beta : I_4 \to I_2$
be the double-point diffeomorphism from a neighbourhood of $p_4$ to a neighbourhood of $p_2$, and
if we let $\alpha_0, \beta_0 \in \Real$ be the coefficient of the linear factor of the Taylor expansions of 
the $x$-coordinates, then the intersection is transverse provided either
$\alpha_0(\beta_0-1) \neq 0$ or $\alpha_0 -1 \neq 0$, moreover the quadruple is the degeneration of 
a $1$-parameter family of downward parabolas if and only if $\alpha_0 \beta_0 \neq 1$.
One convenient way to verify this is to use Lagrange (degree $2$) polynomial interpolation to find
the parabola through three given points, then the equation of the polynomial passing through the fourth
gives the transversality condition.   

There is the degenerate situation where all four points
belong to a common vertical line in $D^2$, they are the degeneration two distinct parabolic quadruples
when $\alpha_0, \beta_0 < 0$ or none (otherwise).  These configurations generically do not occur, but we 
will see for certain spun knots they will show-up on occasion.  Regardless, mod-$2$ they contribute nothing
to the count. There is one further degenerate case of a quadruple
on a degenerate parabola with $p_1$ and $p_4$ infinitesimally separated, as well as $p_2$ and $p_3$, but this
is a high co-dimension phenomena that is avoidable. 

Lastly let's consider neighbourhoods of intersections of the form $\mathcal{C}_{\cone{1}\cthree{2}{3}{4}} \cap 
f_*^{-1}(\mathcal{C}_{\ctwo{1}{3}\ctwo{4}{2}})$, i.e. as limits of the associated parabolic intersections.
First observe all degenerations where the four points collapse to belong to a single vertical line in $D^2$
have co-dimension $\geq 1$ and are therefore avoidable.   We now consider the transversality condition for
intersections of the form in Figure \ref{cd0vt} (b). Let the double-point map denoted by the red 
arrow that sends $p_4$ to $p_2$ be denoted by $\alpha$.  Then the transversality condition for the configuration
in Figure \ref{cd0vt} (b) is given by $\alpha_0 \neq 1$ where $\alpha_0$ is the coefficient of the linear term
in the Taylor expansion of the $x$-coordinate of $\alpha$.  Moreover, these configurations are the degeneration of
a unique $1$-parameter family of parabolas.  Like in the previous case, Lagrange polynomial interpolation 
allows us to reduce these observations to the algebra of quadratic polynomials. 

% Consider perturbation of the x-coordinates of the vertical triple w/ x-coords
%  x, y, \alpha_0 x. The zero of this parabola (away from this triple) goes between max/min {x,y,\alpha_0 x}
% and infinity depending on y, as it goes between x and \alpha_0 x.  That's how I would write down a detailed
% argument, at present. Technically it's not as y goes between x and \alpha_0 x but between x and
% the straight line between (x,altitude) and (\alpha_0 x, altitude). 

\begin{prop}\label{4cyc-form}
Given an element $f \in \Emb(D^2, D^4)$ such that the orthogonal projection 
$\pi : D^4 \to D^3$ to $D^3 \times \{0\}$ is a regular diagram (see discussion following
Figure \ref{parabdiag}) such that the vertical triples and pairs in the domain transversely
intersect the vertical pairs and triples in the co-domain (i.e. the $4$-cycles of overcrossings
and $2$-cycles of overcrossings respectively) then $\mu(f)$ is the mod-$2$ count of the
$4$-cycles and $2$-cycles of overcrossings, respectively.
\end{prop}

One natural extension of the ideas of this section would be to consider $p$-tuples of points
on round circles in $S^j$ being mapped by our embedding $S^j \to S^n$ to $p$-tuples of points
in $S^n$ on a round circle, in  some other cyclic ordering.  We can ensure compactness of the
intersection by similarly demanding that this ordering is not compatible with the cyclic ordering.
For example, if $p$ is prime, we use the standard cyclic ordering in the domain, and the cyclic
ordering given by multiplication by $2$ in $\Zed_p$, in the co-domain.  Thus if we seek a degree
invariant using $k$-parameter families of such embeddings, this would require the equality
$$ j = k + pj - (p-3)(j-1) -(p-3)(n-1)$$
which reduces to
$$ k = (p-3)n-2j-2p+6.$$
Notice in the co-dimension $2$ case this class is an invariant of $\pi_{(p-5)n -2p +10} \Emb(S^j, S^n)$. 
There are often distinct non-consecutive orderings of $\Zed_p$, thus there will often be additional linking
invariants associated to these classes, although we do not compute any here.

\section{Computing $\mu$ on $2$-knots}

In this section we compute $\mu$ on various $2$-knots, including the Yoshikawa table \cite{Yosh}, an infinite family of 
spun knots and an example of Fox.  

Our strategy will be to take a $2$-knot in the form $f : D^2 \to D^4$ and break the computation of $\mu(f)$ into two steps. 
The first step is to compute the double-point diagram of $f$, i.e. the points in $D^2$ where $f$ is two-to-one, i.e. we are essentially
sketching the sets (1), (2), (3) from Section \ref{invtsec} of double and triple points, but represented as a collection of curves 
(and automorphisms of curves) in $D^2$.  These 
are sometimes called fold/decker sets \cite{Cart}. From these diagrams we will compute $\mu$, using
either the formulas from Section \ref{invtsec}, specifically Propositions \ref{halfprop} and \ref{4cyc-form}.

To remind readers, Yoshikawa diagrams of $2$-knots are much like `bridge position' for classical knots.  In bridge position, one has a
linear Morse function on $S^3 \to \Real$, which is Morse also on the knot, moreover all the local maxima are global maxima, similarly
all the local minima are global minima, i.e. occurring at the same altitude. More simply one could say the Morse function
is self-indexing on the knot. For Yoshikawa diagrams we have a linear Morse function
$S^4 \to \Real$ which restricts to a Morse function on the $2$-knot, and all the critical points of index $i$ occur at a common altitude (for each $i = 0, 1, 2$), i.e. the function is self-indexing.  Whereas classical knots in bridge position are described by the braid between the max an min, interestingly
the feature of the Yoshikawa diagram that describes the $2$-knot is the intersection with level $i=1$, plus one small decoration.  
The decoration is a small red dash that indicates how the singularities are resolved as one transitions to nearby level-sets. 
Thus our Yoshikawa diagram is an immersed link in $S^3$ with a number of regular double-points corresponding to the number of 
saddles of the Morse function restricted to the $2$-knot, together with the red decoration of the saddle points.

\begin{figure}[H]
$$\includegraphics[width=3cm]{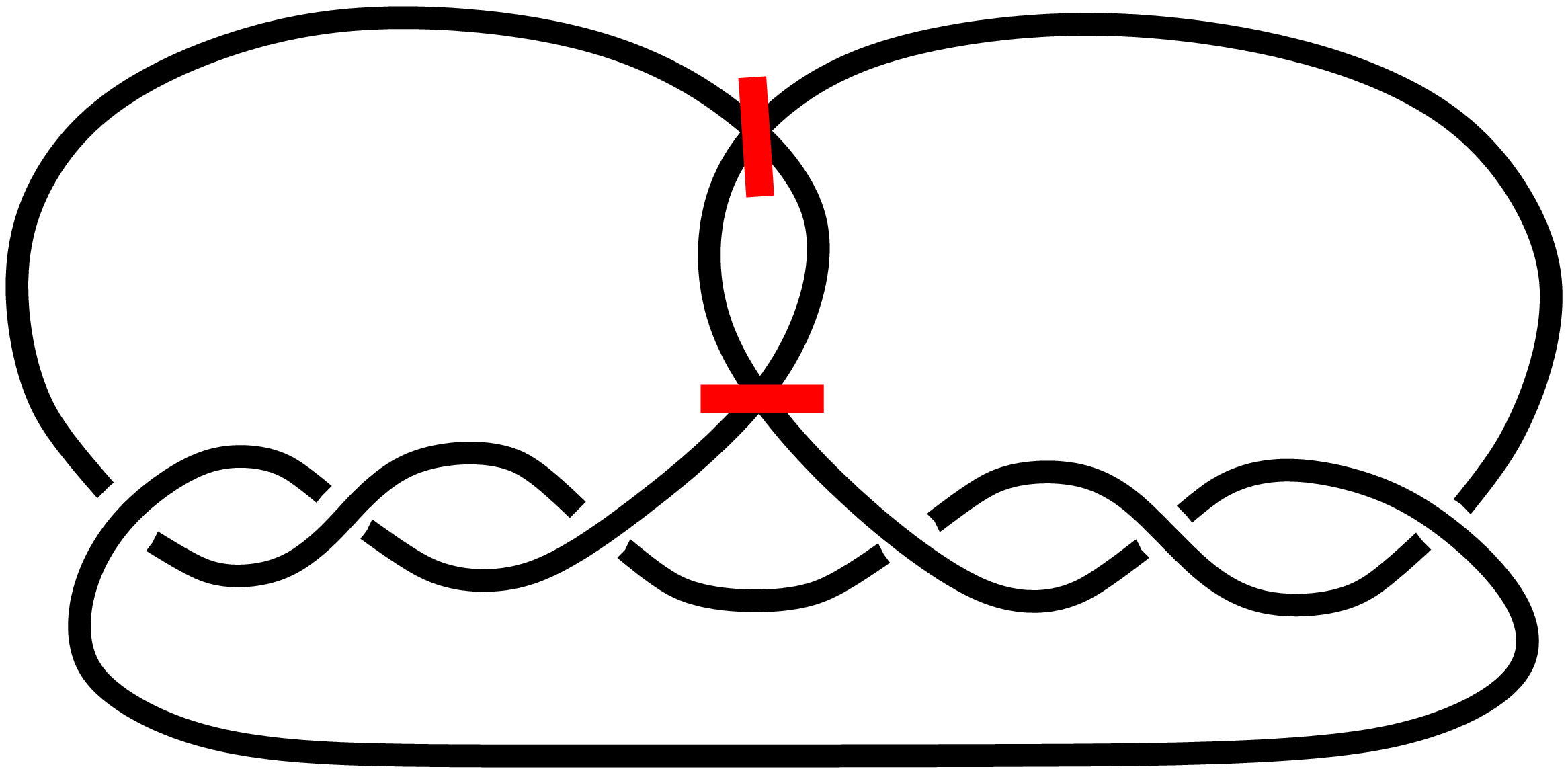}$$ 
\caption{Yoshikawa $8_1$ diagram}\label{81diagram}
\end{figure}

\begin{figure}[H]
{
\psfrag{A}[tl][tl][0.7][0]{(A)}
\psfrag{B}[tl][tl][0.7][0]{(B)}
\psfrag{C}[tl][tl][0.7][0]{(C)}
\psfrag{D}[tl][tl][0.7][0]{(D)}
\psfrag{E}[tl][tl][0.7][0]{(E)}
\psfrag{F}[tl][tl][0.7][0]{(F)}
$$\includegraphics[width=16cm]{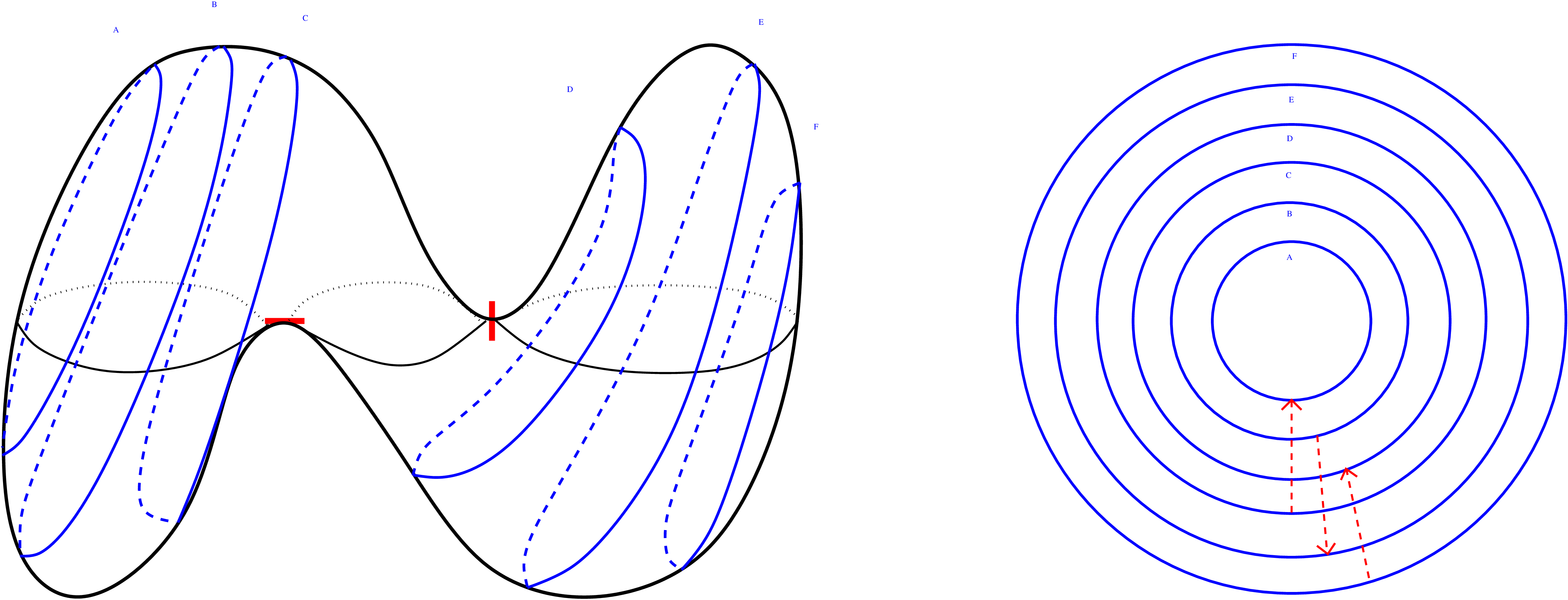}$$
\caption{Double-point diagram for $8_1$ w/Morse height function (left). Projected into plane (right). Red arrows 
depict over-to-under  diffeomorphisms.}\label{dpd81}
}
\end{figure}

If one colours each crossing in Figure \ref{81diagram} blue, and keep track of the crossings as one resolves the Yoshikawa diagram into
$2$-component trivial links, together with the trivialization time parameter represented as up/down, one gets
Figure \ref{dpd81}, left. In this we are only presenting the surface (unknotted) in $\Real^3$ together with the Morse height
function and the double-point sets, i.e. the only detail of the original embedding we are keeping track of is the double-point
set and the height function. In Figure \ref{dpd81} right we have stereographically projected the $S^2$ and identified with
$\Real^2$. The double-point automorphisms are represented by red arrows, and they are radial arrows, with respect to the 
$SO_2$-symmetry in the diagram.

All our crossing diffeomorphisms as diffeomorphisms from the over-to-under crossing curves map (D) to (A), (B) to (E) and (F) to (C). From this we can compute $\mu(8_1)$.  

\begin{prop}\label{81_prop}$\mu(8_1)=0$. 
\begin{proof}
We use the formula for the computation of $\mu$ in Proposition \ref{4cyc-form}.  In this diagram there is the unique 
$4$-cycle of overcrossing at the top of Figure \ref{dpd81} (right), but it is non-transverse.  As described in
the lead-in to Proposition \ref{4cyc-form}, perturbing our detecting manifolds into the space of parabolic quadruples, 
we make the intersection transverse and this will either create a pair of $4$-cycles, or none. In this case, it is none. 

Alternatively, if we use the formula in Proposition \ref{halfprop}, notice we have an entire circle of quadrisecants in
Figure \ref{dpd81}, all non-transverse and radial.  If we apply a small diffeomorphism to the domain of our embedding that perturbs one
of these circles counter-clockwise $\epsilon$-degrees with $\pi > \epsilon > 0$, then the  quadrisecants all vanish, again 
demonstrating $\mu(8_1)=0$. This is the argument used in Corollary \ref{artinspuncor}. Indeed, $8_1$ is an Artin-spun knot, 
so we could have simply cited this Corollary. 
\end{proof}
\end{prop}

We describe a generalization of the Proposition \ref{81_prop} and Corollary \ref{artinspuncor}.  
We begin by reminding readers of the definition of 
Artin-spun \cite{fox}, Twist-spun \cite{Zeeman} and Litherland-spun \cite{Lith} knots.
Artin spinning only takes as input a knot type.  Twist-spinning takes as input a knot type and integer.  
Litherland spinning takes as input a loop in the space of knots, specifically an element of $\pi_0 \Map(S^1, \Emb(D^j, D^n))$. 
One can view these spinning operations as a direct analogue of the braid-closure construction, thought of as a map
$\Omega C_n(\Real^2) \to \Emb(\sqcup_n S^1, \Real^3)$, where $C_n(\Real^2) \equiv \Emb(\{1,2,\cdots,n\}, \Real^2)$ is
the configuration-space of $n$ points in $\Real^2$.  When $\Emb(D^j, D^n)$ is connected, Litherland spinning is known 
to coincide with the connecting map for the 
pseudoisotopy embedding space fibration \cite{BudSurv}, i.e. the natural `graphing' map of the form 
$\Omega \Emb(D^j, D^n) \to \Emb(D^{j+1}, D^{n+1})$.

\begin{figure}[H]
{
\psfrag{v}[tl][tl][0.7][0]{$v$}
\psfrag{w}[tl][tl][0.7][0]{$w$}
\psfrag{fv}[tl][tl][0.7][0]{$f(v)$}
\psfrag{fw}[tl][tl][0.7][0]{$f(w)$}
$$\includegraphics[width=12cm]{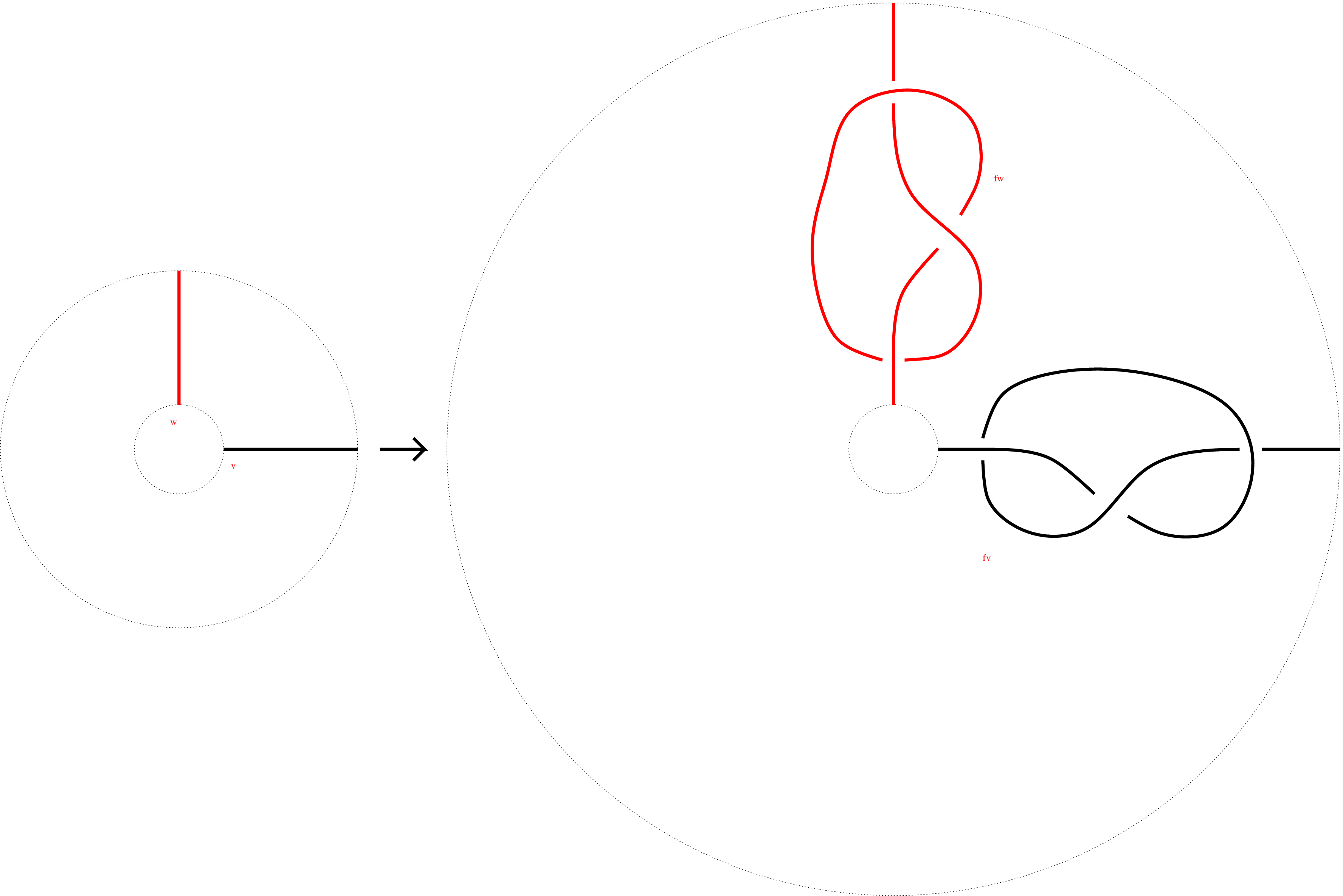}$$ 
}
\caption{A Litherland-spun knot, with $A_v, A_w \in SO_{n+1}$ corresponding to the angles $v,w \in S^1$.}\label{spunknot}
\end{figure}

In the literature typically authors define the exterior 
Litherland-spun knots, rather than defining it at the level of embedding-spaces.  
Let $H^n$ denote a positive half-space in $\Real^n \times \{0\} \subset \Real^{n+1}$. 
We think of $SO_2$ as the subset of $SO_{n+1}$ that fixes $\partial H^n = \Real^{n-1} \times \{0\}^2$ pointwise.
Define the {\bf support} of a map $g : D^j \to D^n$ to be the subset of the domain where $g(p) \neq (p,0)$, i.e. where
$g$ differs from the standard inclusion.  
Given $f \in \Map(S^1, \Emb(D^j, D^n))$ conjugation by an affine-linear map
$\Real^j \to \Real^j$ allows us to assume the support of the embeddings are contained in $H^j \cap D^j$.  
One then defines $F \in \Emb(D^{j+1}, D^{n+1})$ via the formula

$$F(p) = A.\iota f(A, A^{-1}p)$$

where $A \in SO_2$ is the unique element such that $A^{-1}p \in H^j$.  Strictly speaking when $p \in \Real^{j-1}$ 
the term $A$ is ambiguous, but then $F(p)=p$ for any choice of $A$. 
We identify $SO_2$ with $S^1$, which allows us to make sense of the $2^{nd}$ occurence of $A$ in the
expression for $F(p)$. The map $\iota : D^n \to D^{n+1}$  is the standard inclusion $\iota(p)=(p,0)$. 

A Litherland-spun knot we call Artin-spun 
provided $f \in \Map(S^1, \Emb(D^j, D^n))$ is constant in the $S^1$-parameter.  
If $f$ is a $1$-parameter family of rigid motions
applied to a fixed knot, then the Litherland-spun knot $F$ is called Twist-spun.

Proposition \ref{81_prop} and Corollary \ref{artinspuncor} immediately generalizes to the following. 

\begin{thm}$\mu(K)=0$ if $K$ is an Artin-spun $2$-knot.
\begin{proof}
An Artin-spun knot has a double-point diagram similar to Figure \ref{dpd81}, consisting of nested concentric circles with
double-point maps being radial.  Along each radial line the diagram is precisely the chord diagram of the classical knot
we are spinning.  Thus again using either Proposition \ref{4cyc-form} or Proposition \ref{halfprop}, we have $\mu(K)=0$.
\end{proof}
\end{thm}

\begin{figure}[H]
{
\psfrag{x}[tl][tl][0.7][0]{$x$-axis}
\psfrag{y}[tl][tl][0.7][0]{$y$-axis}
$$\includegraphics[width=8cm]{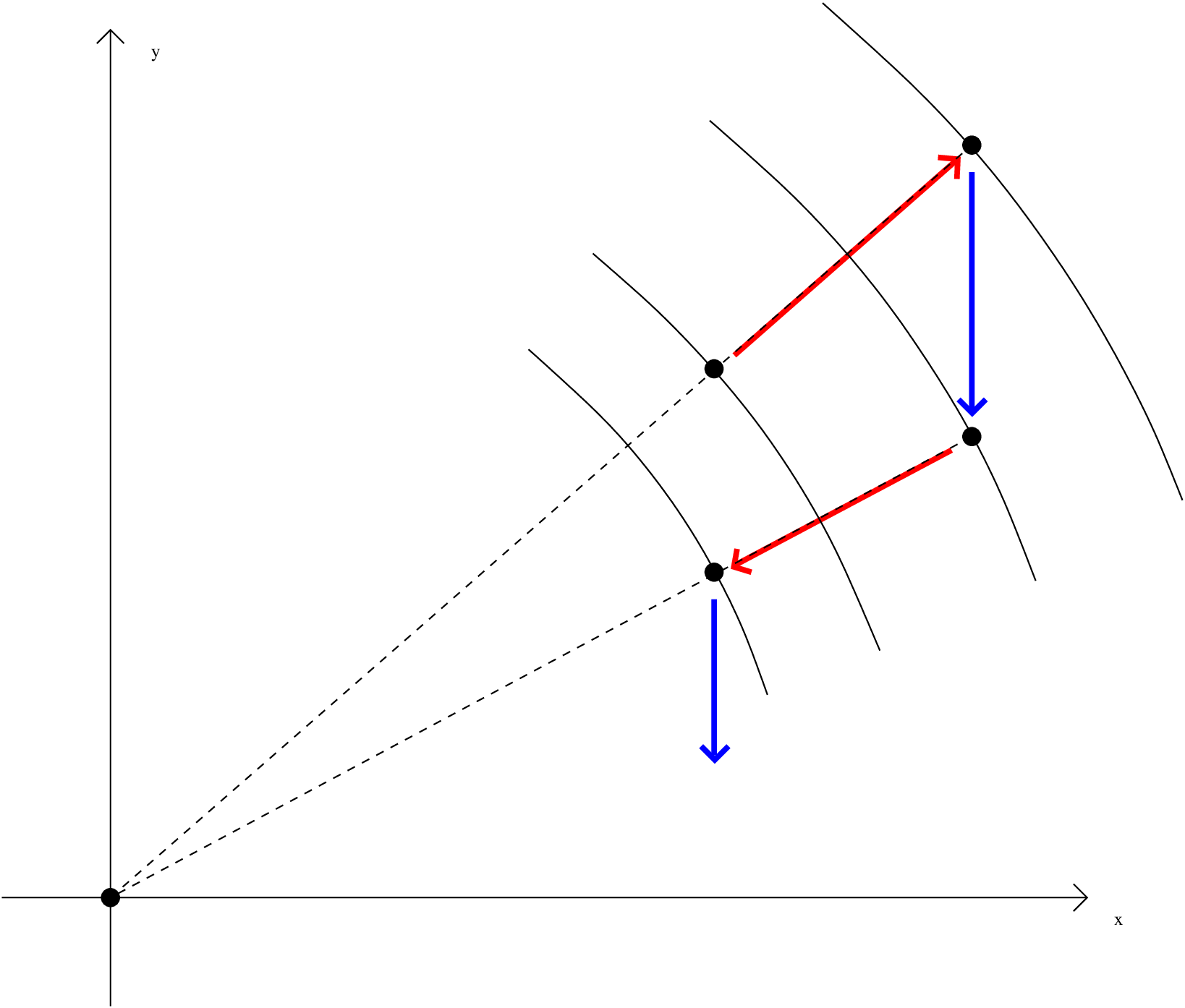}$$ 
}
\caption{The trouble with $4$-cycles of overcrossings for spun knots.}\label{spunquadknot}
\end{figure}

Going one small step further we can prove

\begin{thm}\label{muds}$\mu(K)=0$ if $K$ is Litherland-spun (also known as deform-spun).
\begin{proof}
For this argument we use the formula in Proposition \ref{4cyc-form}.  As before, this diagram consists of $1$-parameter
family of chord diagrams for the knots of $f \in \Map(S^1, \Emb(D^1, D^3))$, radially about the origin.
The only $4$-cycles or $2$-cycles of overcrossings are therefore on the $y$-axis, see for example Figure \ref{spunquadknot}.
These are non-transverse and resolve to an even number of $4$ and $2$-cycles.
\end{proof}
\end{thm}

\begin{figure}[H]
$$\includegraphics[width=4cm]{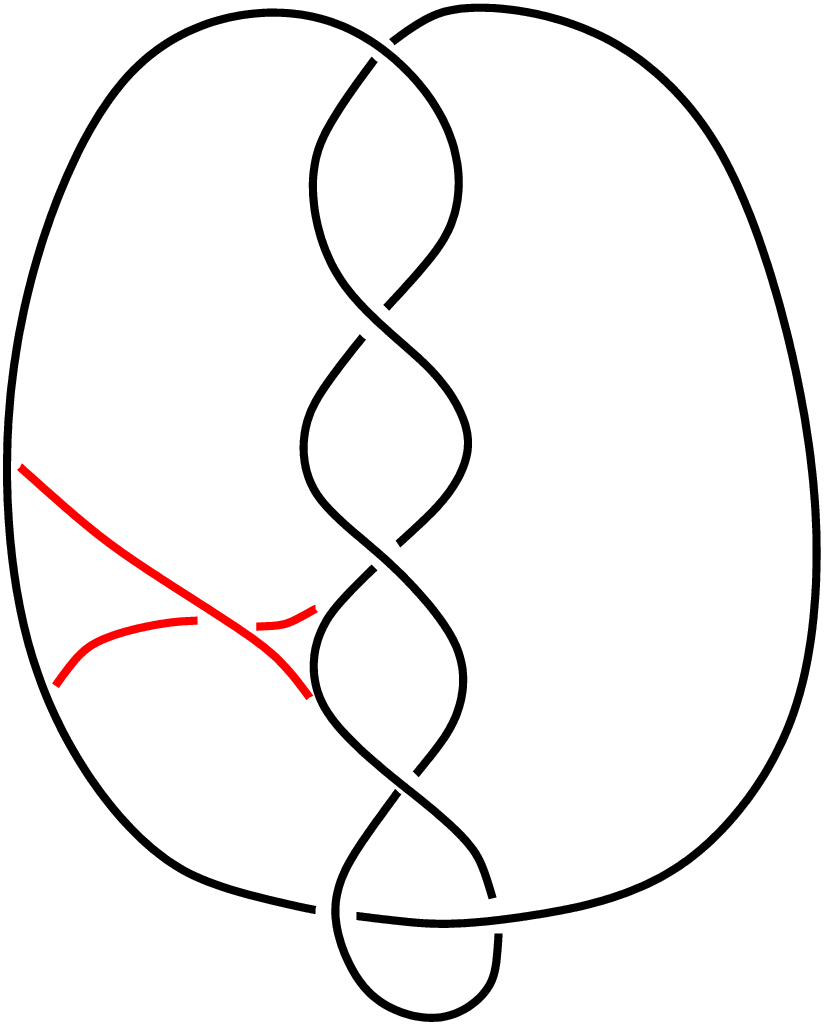}$$ 
\caption{Fox's Quick Trip example 10.}\label{qt10eg}
\end{figure}

Theorem \ref{muds} points us towards non-trivial examples.  For example in \cite{BM} it's noted that a
deform-spun $2$-knot has a symmetric Alexander polynomial.  Thus if we want $\mu(K)=1$ we need a knot with
a non-symmetric Alexander polynomial.  Fox's Quick Trip \cite{fox} Example 10 is one of the simplest $2$-knots
with non-symmetric Alexander polynomials.

Figure \ref{qt10eg} depicts the Stevedore knot in black.  In red one sees the outline of a $1$-handle attachment, 
which from this perspective includes a half-twist.  If one performs embedded
surgery along this $1$-handle it converts the Stevedore knot into a $2$-component trivial link.  
Thus the diagram is depicting a position for a smooth slice disc $D^2 \to D^4$ for the Stevedore knot, 
with the $1$-handle corresponding to the sole critical point of index $(1,1)$ with respect to the radial distance
height function. The height function has two minima (non-degenerate), which are the only other critical points. 
The maximum ($r=1$) is the boundary Stevedore knot. 
The double of this smooth slice disc is Fox's Quick Trip Example 10.  

\begin{figure}[H]
{
\psfrag{A}[tl][tl][0.7][0]{A}
\psfrag{B}[tl][tl][0.7][0]{B}
\psfrag{C}[tl][tl][0.7][0]{C}
\psfrag{D}[tl][tl][0.7][0]{D}
\psfrag{E}[tl][tl][0.7][0]{E}
\psfrag{F}[tl][tl][0.7][0]{F}
\psfrag{G}[tl][tl][0.7][0]{G}
$$\includegraphics[width=12cm]{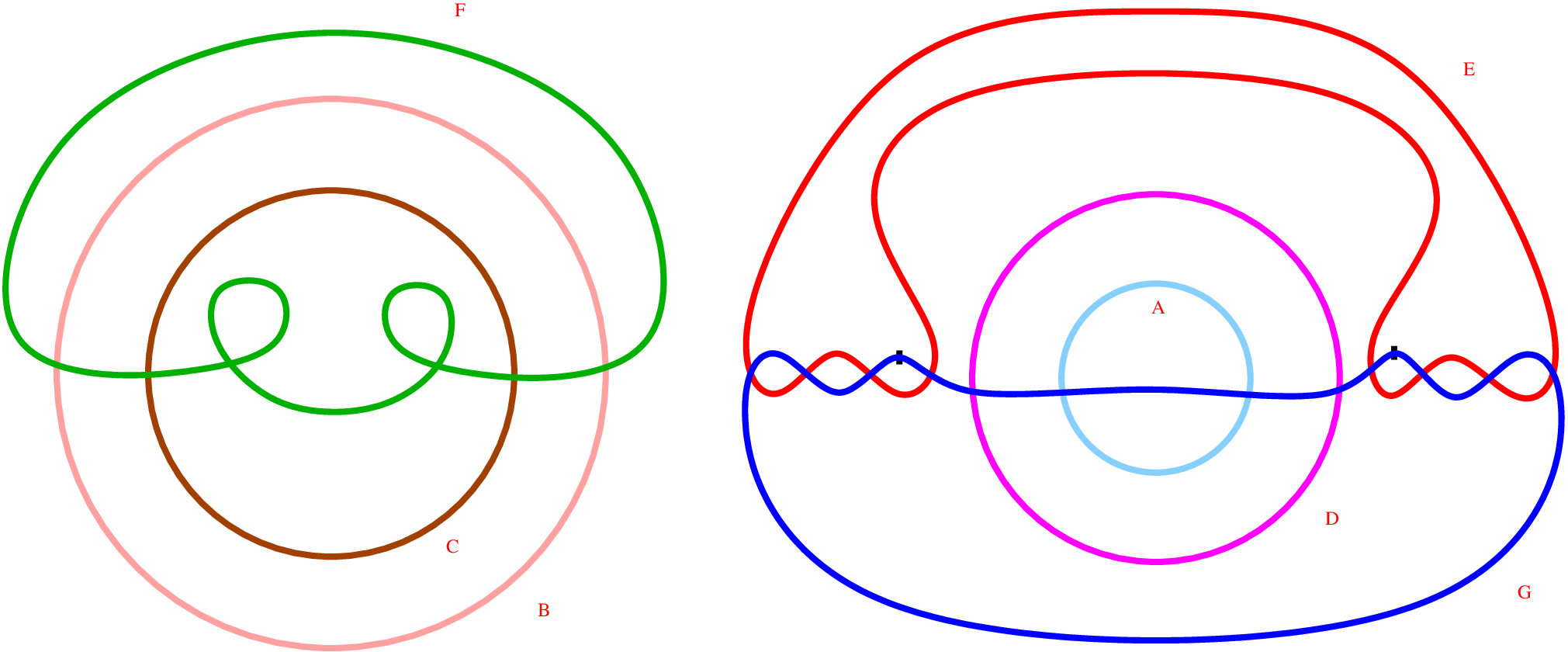}$$ 
}
\caption{Fox's Quick Trip example 10, double-point diagram}\label{qt10eg.1}
\end{figure}

We interpret this knot as an element of $\Emb(D^2, D^4)$, and its double-point diagram in $D^2$ is given by Figure \ref{qt10eg.1}.  
To see the computation, observe that the double-point 
curve $E$ is over $F$, we denote this by $E \to F$.   Similarly, 
$C \to D$ and $A \to B$. The double-point curve $G$ is over itself, i.e. the overcrossing relation is an involution of
the curve $G$, with fixed points marked in black.  We claim there is a unique $4$-cycle of overcrossings in the diagram, 
and no $2$-cycles of overcrossings.  The $4$-cycle has the form $A \textcolor{red}{\to} B \textcolor{blue}{\to} C \textcolor{red}{\to} D \textcolor{blue}{\to} A$. 
%$\xymatrix{A \ar[r] & B \ar[r] & C \ar[r] & D \ar@/_1.0pc/@[red][lll]}$

\begin{prop}If $\text{Fox}_{10}$ denotes the Fox Quick Trip Example 10, then 
$$\mu(\text{Fox}_{10}) = 1.$$
\end{prop}

Our formula for $\mu$ indicates the potential for there to be an analogue of Skein relations for invariants
of knotted $2$-spheres and surfaces in $S^4$.  The idea of Skein-type invariants was explored in 1982 by Giller \cite{CG}.  
An interesting observation
of Giller's is that if one takes a double-point diagram for a knotted $2$-sphere in $\Real^4$, and if one changes one over/under
crossing curve (from over to under or under to over), this is not always the diagram for a $2$-sphere in $\Real^4$.  We would
of course prefer `Skein relation' for knotted surface diagrams in $\Real^4$ to not involve realizability questions.  One of the
appeals of Skein relations is the observation that one can monotonically change crossings to turn any $1$-knot into the
unknot, thus whatever `Skein relation' should mean for $2$-knots, it should not be cast entirely in the language of crossing 
changes and double-point resolutions. 

\section{The remaining cases: both $n > 4$ and $j > 1$}

Arone and Turchin \cite{AT} studied the closely-related space $\overline{\Emb}(D^j, D^n)$, this space is the homotopy-fibre of
the Smale-Hirsch map $\Emb(D^j, D^n) \to \Omega^j V_{n,j}$.  In that paper they describe the rational homotopy groups of
$\overline{\Emb}(D^j, D^n)$ when $n > j + 2$ and most relevant to this paper they show that

$$\Rat \otimes \pi_{2(n-j-2)} \overline{\Emb}(D^j, D^n) \simeq 
\begin{cases} \Rat & \text{if } j=1 \text{ or both } n \text{ and } j \text{ odd}\\ 
              0 & \text{otherwise.}\end{cases}.$$
              
The copy of $\Rat$ above we will simply call the Arone-Turchin class. Arone and Turchin compute the rank of 
$\Rat \otimes \pi_{2(n-j-2)} \Emb(D^j, D^n)$, in particular they notice
that $\Rat \otimes \pi_{2(n-j-2)} \overline{\Emb}(D^j, D^n)$ injects faithfully into $\Rat \otimes \pi_{2(n-j-2)} \Emb(D^j, D^n)$
for all $n \geq j+3$.  

\begin{conj}
The image of the Arone-Turchin class in $\pi_{2(n-j-2)} \Emb(D^j, D^n)$ is detected by the invariant $\mu$ for all $n > j+2$ 
with both $j$ and $n$ odd. 
\end{conj}

The extension problem for the fibration $\Emb(D^j, D^n) \to \Omega^j V_{n,j}$ is subtle and not much is known about it at present.
In the Arone-Turchin paper \cite{AT} they solve the problem rationally, relying on the fact that $V_{n,j}$ is rationally a fairly small 
space.  It should be noted that only a few examples of homotopically non-trivial Smale-Hirsch maps $\Emb(D^j, D^n) \to \Omega^j V_{n,j}$ are known. A recent result on this topic is the paper of Crowley, Schick and Steimle \cite{CSS} 
where they  show this map is non-trivial for $n=j=11$ on the $5$-th homotopy group, i.e. $\pi_5 \Diff(D^{11}) \to \pi_{16} O_{11}$.
In co-dimension $2$ it's known that the immersions $D^{n-2} \to D^n$ realizable as embeddings are precisely those whose $J$-invariant
($J : \pi_{n-2} SO_n \to \pi_n S^n$) is zero. Moreover, this only occurs when $n$ is congruent to $1$ mod $4$ \cite{HM}.

\begin{conj}
The invariant 
$$\mu :\pi_{2(n-j-2)} \Emb(D^j, D^n) \to 
\begin{cases} \Zed & \text{if } j=1 \text{ or both } n \text{ and } j \text{ odd}\\ 
              \Zed_2 & \text{otherwise.}\end{cases}$$
is an epimorphism for all $n \geq j+2$ with $j \geq 1$.
\end{conj}

As we have noted, this conjecture is known to be true for all $n \geq 3$ with $j=1$, as well as the case $(n,j)=(4,2)$. 
A potential starting point to resolve this conjecture would be the cycles constructed by Sakai and Watanabe \cite{SW}.

\begin{conj}
When $n=j+2$ the invariant
$$\mu : \pi_{0} \Emb(D^j, D^n) \to 
\begin{cases} \Zed & \text{if } j=1 \text{ or both } n \text{ and } j \text{ odd}\\ 
              \Zed_2 & \text{otherwise.}\end{cases}$$
can be expressed in terms of the rational Alexander modules. 
\end{conj}

This conjecture is asserting that $\mu$ is computable in terms of the rational homology of the
universal abelian cover of the knot exterior, as a module over $\Rat[t^\pm]$.  At present there is limited
evidence for this, but for the handful of knots where $\mu$ has been computed and the Alexander module is
$\Zed$-torsion, $\mu$ appears to be zero.  One natural conjecture would be that $\mu$ is a (sometimes mod-$2$)
rational Alexander module Euler characteristic, as in \cite{BCSS}, \cite{PV98}.

Habiro, Kanenobu and Shima \cite{HKS} have a notion of rational finite-type invariant for ribbon $2$-knots and prove that
finite-type invariants are polynomial functions in the coefficients of the Alexander polynomial.  Likely our result 
should fit into a broader theory of finite-type invariants, but at present ours is torsion valued and defined for
all knots, not just ribbon 2-knots.  This author is currently unaware of a satisfying definition of finite-type invariant
for arbitrary $2$-knots, i.e. there is no full analogue to the works of Birman-Lin \cite{BL} or Vassiliev \cite{V} 
for co-dimension two knot theory above dimension $3$. 

Recently Gauniyal and Turchin \cite{GT} have used similar techniques as in this paper to compute invariants of Haefliger knots, 
i.e. such as the isomorphism $\pi_0 \Emb(S^3, S^6) \to \Zed$. Like us, they use the cobordism class of the double-point set.  
Unlike us they keep their invariants largely in the language of those cobordism classes, constructing an invariant of
$\pi_0 \Emb(S^j, S^n)$, for $n-j>2$, whereas we `double down' on double-point formulas, which forces our invariant to
live in the homotopy group $\pi_{2(n-j-2)} \Emb(S^j, S^n)$.   Both our works are guided by Arone and Turchin's rational
homotopy computations \cite{AT}, but ours is only loosely guided in that we happily produce a torsion invariant.  

\providecommand{\bysame}{\leavevmode\hbox to3em{\hrulefill}\thinspace}

%\Addresses

\end{document}